\documentclass{amsart}
\input{xypic.tex} 
\xyoption{all}
\usepackage{ifthen,amsthm,amssymb,latexsym,amsmath}
\usepackage[all]{xy}
\newcommand{\cala}{{\mathcal A}}
\newcommand{\calv}{{\mathcal V}}
\newcommand{\calp}{{\mathcal P}}
\newcommand{\oy}{{\mathcal O}_Y}

\newcommand{\id}{{\mathbf 1}}

\newcommand{\R}{\mathbb{R}} \newcommand{\C}{\mathbb{C}} \newcommand{\F}{\mathbb{F}}
\newcommand{\Z}{\mathbb{Z}}  \newcommand{\A}{\mathbb{A}}
\newcommand{\p}[1]{{\mathbb{P}^{#1}}} \newcommand{\op}{{\mathcal O}_{\mathbb{P}^{2}}}
\newcommand{\oo}{{\mathcal O}}

\DeclareMathOperator{\codim}{{codim}}
\DeclareMathOperator{\rk}{{rank}}
\DeclareMathOperator{\im}{{im}}

\newtheorem{coro}{Corollary}[section]
\newtheorem{defi}[coro]{Definition}

\newtheorem{lemm}[coro]{Lemma}
\newtheorem{prop}[coro]{Proposition}
\newtheorem{rema}[coro]{Remark}
\newtheorem{theo}[coro]{Theorem}

\newcommand{\marginnote}[1]{\ifthenelse{\isodd{\thepage}}{\normalmarginpar}
{\reversemarginpar}\marginpar{\fbox{\parbox{24mm}{\sloppy\footnotesize #1}}}}

\begin{document}

\title{The ADHM Variety and Perverse Coherent Sheaves}

\author{Marcos Jardim}
\author{Renato Vidal Martins}

\address{IMECC - UNICAMP \\
Departamento de Matem\'atica \\ Caixa Postal 6065 \\
13083-970 Campinas-SP, Brazil}
\email{jardim@ime.unicamp.br}
\address{ICEx - UFMG \\
Departamento de Matem\'atica \\ 
Av. Ant\^onio Carlos 6627 \\
30123-970 Belo Horizonte MG, Brazil}
\email{renato@mat.ufmg.br}

\begin{abstract}
We study the full set of solutions to the ADHM equation as an affine algebraic set, the ADHM variety. We gather the points of the ADHM variety into subvarieties according to the dimension of the stabilizing subspace. We compute dimension, and analyze singularity and reducibility of all of these varieties. We also connect representations of the ADHM quiver to coherent perverse sheaves on $\mathbb{P}^2$ in the sense of Kashiwara.
\end{abstract}

\maketitle


\section{Introduction}

Atiyah, Hitchin, Drinfeld and Manin showed in \cite{ADHM}, using Ward correspondence and algebro-geometric techniques (monads) introduced by Horrocks, that all self-dual connections on euclidean 4-dimensional space (a.k.a. instantons) have a unique description in terms of linear algebra. 

A few years later in \cite{D1}, Donaldson restated the ADHM description in terms of the following data. Let $V$ and $W$ be complex vector spaces, with dimensions $c$ and $r$, respectively. Let $A,\,B\in{\rm Hom}(V,V)$, let $I\in {\rm Hom}(W,V)$ and let $J\in {\rm Hom}(V,W)$. Consider the following equations:
\begin{align}
[A,B] + IJ & =  0 \label{adhm1} \\
[A,A^\dagger]+[B,B^\dagger] + II^\dagger - J^\dagger J &=  0. \label{adhm2}
\end{align}
The group $GL(c)$ acts on the set of solutions of the first equation, sending a datum $(A,B,I,J)$ to $(gAg^{-1},gBg^{-1},gI,Jg^{-1})$, where $g\in GL(c)$, while the unitary group $U(c)$ preserves the second equation. Among other things, Donaldson showed that the \emph{regular} (Definition \ref{m-def}) solutions to equation (\ref{adhm1}) modulo the action of $GL(c)$ parametrize the moduli space of holomorphic bundles on $\C\p2$ that are framed at a line, so-called line at infinity. The key observation is that the set of regular solutions of the first equation modulo $GL(c)$ is isomorphic to the set of regular solutions of both equations modulo $U(c)$, which in turn is identified with the moduli space of framed instantons on $\R^4$.

More recently, Nakajima showed that relaxing the regularity condition to a weaker stability condition one obtains the  moduli space of framed torsion free sheaves on $\C\p2$, see \cite{N}. In particular, Nakajima also obtained a linear algebraic description of the Hilbert scheme of points on $\C^2$ and showed that it admits a hyperk\"ahler structure.

Generalizations of the ADHM construction have given rise to a wide variety of important results in many different areas of mathematics and mathematical physics. For instance, Kronheimer and Nakajima constructed instantons on the so-called ALE manifolds \cite{KN}. Nakajima then proceed to further generalize the construction, and studied quiver varieties and representations of Kac-Moody Lie algebras \cite{N1}. More recently, the ADHM construction of instantons was adapted to construct analogues of Yang-Mills instantons in noncommutative geometry \cite{NS} and supergeometry \cite{T}, while the ADHM construction of holomorphic bundles on $\C\p2$ was generalized to the construction of instantons bundles on $\C\mathbb{P}^n$ \cite{FJ2,J}. 

We start this paper by studying the full set of solutions of (\ref{adhm1}), called the \emph{ADHM equation}, as an affine algebraic variety $\calv(r,c)$, called the \emph{ADHM variety}. Our first concern are the points $X\in\calv(r,c)$ which are not {\it stable}. In order to study them, we define the {\it stabilizing subspace} $\Sigma_X$ as the subspace of $V$ to which the restriction of $X$ is stable (c.f. Definition \ref{defsig}). We break down the ADHM variety into disjoint subsets $\calv(r,c)^{(s)}:=\{X\,|\,\dim\Sigma_X=s\}$; in this sense, the set of stable points $\calv(r,c)^{\rm st}$ now corresponds to $\calv(r,c)^{(c)}$. It turns out that the study of these subsets happens to give some relevant properties of the whole variety itself. The results we get are:

\

\noindent {\bf Theorem 1.} \emph{The ADHM variety $\calv(r,c)$ is a set theoretic complete intersection which is irreducible if and only if $r\geq 2$. Moreover, the set $\calv(r,c)^{(s)}$ is an irreducible quasi-affine variety of dimension $2rc+c^2-(r-1)(c-s)$ which is nonsingular if and only if either $s=c$ or $s=c-1$. In particular, $\calv(r,c)^{\rm st}$ is a nonsingular irreducible quasi-affine variety of dimension $2rc+c^2$.}

\

The reader will note that everything done concerning stable points could have been said for {\it costable} as well. In fact, part of this -- {\it costabilizing subspaces}, for instance -- actually appears here in the proof of Proposition \ref{prpule}, which provides an analogue of the so-called \emph{Gieseker to Uhlenbeck map} (c.f. \cite[pp. 196-198]{B}) purely in terms of the ADHM data.

For the second part of the paper, it is important to put the subject within a categorical framework. We pass from ``points of the ADHM variety $\calv$'' to ``objects in the category $\mathcal{A}$ of representations of the \emph{ADHM quiver}''. Such objects are of the form $\mathbf{R}=(V,W,X)$ where $V$, $W$ and $X$ are as above. For every $\mathbf{R}\in\mathcal{A}$ we define a \emph{stable restriction} by $\mathbf{S}_{\mathbf{R}}:=(\Sigma_X,W,X|_{\Sigma_X})$ and a \emph{quotient representation} by $\mathbf{Z}_{\mathbf{R}}:=(V/\Sigma_X,\{0\},(A',B',0,0))$ with $A',B'$ commuting operators in $V/\Sigma_X$, see (\ref{equzzr}). The triple $(\dim W,\dim \Sigma_X,\codim \Sigma_X)$ is called the \emph{type vector} of $\mathbf{R}$.

From every $\mathbf{R}\in\mathcal{A}$ one constructs a complex $E_{\mathbf{R}}^{\bullet}\in \text{Kom}(\mathbb{P}^2)$ (Definition \ref{defcom}), called the \emph{ADHM complex} associated to $\mathbf{R}$. Nakajima has shown in \cite[Chp.  2]{N}, based on previous constructions due to Barth and Donaldson \cite{D1}, that if $\mathbf{R}$ is stable, then $E:=H^0(E^\bullet_\mathbf{R})$ is the only nontrivial cohomology sheaf of the ADHM complex; moreover, $E$ is a torsion free sheaf such that $E|_{\ell}\simeq W\otimes\oo_{\ell}$ for a fixed line $\ell\subset\mathbb{P}^2$. Conversely, given any torsion free sheaf $E$ on $\p2$ whose restriction to $\ell$ is trivial, then there is a stable solution of the ADHM equation such that $E$ is isomorphic to the cohomology sheaf of the corresponding ADHM complex. Besides, $H^0(E^\bullet_\mathbf{R})$ is locally free if and only if $\mathbf{R}$ is regular.

Following ideas of Drinfeld, Braverman, Finkelberg and Gaitsgory  introduced in \cite[Sec. 5]{BFG} the notion of perverse coherent sheaves, and showed, generalizing \cite[Chp. 2]{N}, that such objects correspond to arbitrary solutions of the ADHM equation. Here, we show that the complex associated to such solutions can also be considered perverse in the sense of Kashiwara's ``family of supports" approach. More precisely, we introduce, following \cite{Ka}, a perverse $t$-structure on $D^{\rm b}(\p2)$; objects in the core of such $t$-structure, denoted $\calp$, are what we call \emph{coherent perverse sheaves on $\mathbb{P}^2$}. We characterize objects $E^\bullet$ in $\calp$, and define the \emph{rank, charge} and \emph{length} of such an $E^\bullet$ as, respectively, $\rk (H^0(E^\bullet))$, $c_2(H^0(E^\bullet))$ and $\text{length}(H^1(E^\bullet))$. We then prove:

\

\noindent {\bf Theorem 2.} \emph{If $\mathbf{R}$ is a representation of the ADHM quiver with type vector $(r,s,l)$, then the associated ADHM complex $E^\bullet_\mathbf{R}$ is a perverse coherent sheaf on $\p2$ of rank $r$, charge $s$ and length $l$.
Moreover, 
\begin{itemize}
\item[(i)] $H^0(E^\bullet_{\mathbf{R}})\simeq H^0(E^\bullet_{\mathbf{S}_{\mathbf{R}}})$;
\item[(ii)] $H^1(E^\bullet_{\mathbf{R}})\simeq H^1(E^\bullet_{\mathbf{Z}_{\mathbf{R}}})$;
\item[(iii)] $H^0(E^\bullet_{\mathbf{R}})$ is locally free if and only if $\mathbf{R}$ is costable.
\end{itemize}}

Recently, Hauzer and Langer \cite{HL} solved this problem for the complex ADHM equations introduced in \cite{FJ2}, showing that arbitrary solutions of the complex ADHM equation correspond to perverse instanton sheaves on $\p3$. This picture can the further generalized, and, as it was pointed out in \cite{J}, an analogous relation between arbitrary solutions of the $d$-dimensional ADHM equations and perverse instanton sheaves on $\mathbb{P}^{d+2}$ should also hold. 

\

\noindent{\bf Acknowledgments.} The first named author is partially supported by the CNPq grant number 305464/2007-8 and the FAPESP grant number 2005/04558-0.


\section{The ADHM Data}

Let $V$ and $W$ be complex vector spaces, with dimensions $c$ and $r$, respectively. The {\em ADHM data} is the set (or space) given by
$$ 
\mathbf{B}=\mathbf{B}(r,c) := {\rm Hom}(V,V)\oplus{\rm Hom}(V,V)\oplus {\rm Hom}(W,V)\oplus{\rm Hom}(V,W).
$$
An element (or point) of $\mathbf{B}$ is a {\em datum} $X=(A,B,I,J)$ with $A,B\in{\rm End}(V)$, $I\in{\rm Hom}(W,V)$ and $J\in{\rm Hom}(V,W)$.

\begin{defi} 
\label{m-def}
A datum $X=(A,B,I,J)\in\mathbf{B}$ is said
\begin{enumerate}
\item[(i)] {\em stable} if there is no subspace $S\subsetneqq V$ with $A(S),B(S),I(W)\subset S$;
\item[(ii)] {\em costable} if there is no subspace $0\neq S\subset V$ with $A(S),B(S)\subset S\subset \ker J$;
\item[(iii)] {\em regular} if it is both stable and costable.
\end{enumerate}
We call $\mathbf{B}^{\rm st}$, $\mathbf{B}^{\rm cs}$ and $\mathbf{B}^{\rm reg}$ the sets of stable, costable and regular data, respectively.
\end{defi}

If $T$ is a linear map, let $T^\dagger$ be its adjoint operator. For $X=(A,B,I,J)\in\mathbf{B}$ one defines
$$ X^\star := (B^\dagger,-A^\dagger,J^\dagger,-I^\dagger). $$
It is easily seen that $(X^\star)^\star=-X$, and that $X$ is stable if and only if $X^\star$ is costable. The anti-linear endomorphism of $\mathbf{B}$ given by $X\to X^\star$ plays the role of the duality needed in this work, something that will be more clear later on in Section \ref{secprv}.

We are specially interested in the following morphism
\begin{gather*}
\begin{matrix}
\mu : & \mathbf{B} & \longrightarrow &{\rm End}(V)   \\
      & (A,B,I,J)  & \longmapsto     & [A,B] + IJ
\end{matrix}
\end{gather*}
which will define the variety we will deal with in next section. It is easily checked that, for any $X=(A,B,I,J)\in\mathbf{B}$, the derivative is
\begin{gather*}
\begin{matrix}
D_X\mu : & \mathbf{B} & \longrightarrow &{\rm End}(V)   \\
      & (a,b,i,j)  & \longmapsto     & [A,b] + [a,B] + Ij + iJ.
\end{matrix}
\end{gather*}

\begin{lemm}
\label{lemsj1}
Let $X=(A,B,I,J)\in\mathbf{B}$. Then $D_X\mu$ is surjective iff the map
\begin{gather*}
\begin{matrix}
{\rm End}(V) & \longrightarrow & \mathbf{B}   \\
y & \longmapsto     & ([A,y],[B,y],yI,Jy)
\end{matrix}
\end{gather*}
is injective.
\end{lemm}

\begin{proof}
The map $D_X\mu$ is not surjective iff there is a nonzero $y\in{\rm End}(V)$ such that $y^\dagger\in(\im D_X\mu)^\perp$. But $y^\dagger\in(\im D_X\mu)^\perp$ iff, for every $(a,b,i,j)\in\mathbf{B}$, holds
$$ 
{\rm Tr}(D_X\mu(a,b,i,j)y)=0
$$
which is equivalent, for every $(a,b,i,j)\in\mathbf{B}$, to the following equalities
\begin{align*}
{\rm Tr}([A,b]y)={\rm Tr}([y,A]b)&=0 \\
{\rm Tr}([a,B]y)={\rm Tr}(a[B,y])&=0 \\
     {\rm Tr}((Ij)y)={\rm Tr}((yI)j)&=0 \\
     {\rm Tr}((iJ)y)={\rm Tr}(i(Jy))&=0 
\end{align*}
which hold if and only if 
$$
[A,y]=[B,y]=yI=Jy=0
$$
and we are done.
\end{proof}

\begin{prop}
\label{prpsrj}
If $X$ is stable or costable then $D_{X}\mu$ is surjective.
\end{prop}

\begin{proof}
Set $X=(A,B,I,J)$. If $D_{X}\mu$ is not surjective, from Lemma \ref{lemsj1} there is a nonzero $y\in{\rm End}(V)$ such that $[A,y]=[B,y]=yI=Jy=0$. 

Since $yI=0$, then $I(W)\subset\ker y$, while $[A,y]=[B,y]=0$ implies that $\ker y$ is $A$- and $B$-invariant. But $\ker y$ is a proper subspace of $V$ because $y\neq 0$. Hence $X$ is not stable.

Similarly, since $Jy=0$, then $\im y\subset \ker J$, while $[A,y]=[B,y]=0$ implies that $\im y$ is $A$- and $B$-invariant. But $\im y$ is a nonzero subspace of $V$ because $y\neq 0$. Hence $X$ is not costable as well.
\end{proof}

The group $G:=GL(V)$ acts naturally on $\mathbf{B}$. For $g\in G$ and $X=(A,B,I,J)\in\mathbf{B}$ one defines
$$
g\cdot X := (gAg^{-1},gBg^{-1},gI,Jg^{-1}).
$$
One also defines the stabilizer subgroup of $X$ as
$$
G_X:=\{g\in G\,|\,g\cdot X=X\}.
$$

\begin{prop}
\label{trstab}
If $D_{X}\mu$ is surjective then $G_X$ is trivial.
\end{prop}

\begin{proof}
Let $X=(A,B,I,J)\in\mathbf{B}$. If $G_X$ is nontrivial, let $g\ne\id$ be such that $gA=Ag$, $gB=Bg$, $gI=I$ and $Jg=J$. Hence
$$
[A,g-\id]=[B,g-\id]=(g-\id)I=J(g-\id)=0.
$$
Since $g-\id\neq 0$, Lemma \ref{lemsj1} implies that $D_X\mu$ is not surjective.
\end{proof}

Propositions \ref{prpsrj} and \ref{trstab} led us considering the sets 
$$ 
\mathbf{B}^{\rm sj} = \{ X\in\mathbf{B}\,|\, D_X\mu {\rm ~is~surjective} \}
$$
$$
\mathbf{B}^{\rm ts} = \{ X\in\mathbf{B}\,|\, G_X\ {\rm is~trivial} \} 
$$
which satisfy the sequence of inclusions
$$ 
(\mathbf{B}^{\rm st}\cup\mathbf{B}^{\rm cs}) \subset \mathbf{B}^{\rm sj} \subset \mathbf{B}^{\rm ts}.
$$
which, in general, are strict ones (see Rem. \ref{remstr} below).

\begin{defi}
\label{defsig}
The \emph{stabilizing subspace} $\Sigma_X$ of an ADHM datum $X=(A,B,I,J)$ is the intersection of all subspaces $S\subseteq V$ such that $A(S),B(S),I(W)\subset S$.
\end{defi}

If $S\subseteq V$ satisfies $A(S),B(S),I(W)\subset S$, then one may consider 
$$
X|_S:=(A|_S,B|_S,I,J|_S)\in{\rm End}(S)^{\oplus 2}\oplus{\rm Hom}(W,S)\oplus{\rm Hom}(S,W).
$$
It is clear that $X|_{\Sigma_X}$ is stable and this justifies the term we use. 

For $0\le s\le c$, we define the sets
$$ 
\mathbf{B}^{(s)} := \{ X\in\mathbf{B} ~|~ \dim \Sigma_X= s \}
$$
$$ 
\mathbf{B}^{[s]} := \{ X\in\mathbf{B} ~|~ \dim \Sigma_X\leq s \}
$$
which, respectively, give a disjoint decomposition and filtration of $\mathbf{B}$. Clearly, $\mathbf{B}^{(c)}$ coincides with the set of stable data, while $\mathbf{B}^{(0)}$ is the set of data $X$ for which $I=0$.

Given an ADHM datum $X=(A,B,I,J)$, consider the map 
\begin{gather*}
\begin{matrix}
R(X) : & W^{\oplus c^2}                              & \longrightarrow &  V  \\
      & \bigoplus\limits_{0\leq k,l\leq c-1} w_{kl}  & \longmapsto     & \sum\limits_{0\leq k,l\leq c-1} A^kB^lIw_{kl}.
\end{matrix}
\end{gather*}
It is easily seen that the assignment $R:\mathbf{B}\to {\rm Hom}(W^{\oplus c^2},V)$ defines a regular map. 

\begin{prop}
\label{l1}
For any $X\in\mathbf{B}$ hold
\begin{enumerate}
\item[(i)] $\Sigma_X\supseteq\im R(X)$; 
\item[(ii)] if $R(X)$ is surjective, then $X$ is stable.
\end{enumerate}
\end{prop}

\begin{proof}
If $S\subseteq V$ satisfies $A(S),B(S),I(W)\subset S$, then $\im R(X)\subseteq S$ so (i) follows. Now if $R(X)$ is surjective then $c={\rm rk}\,R(X)\leq\dim \Sigma_X\leq c$. Thus $\dim\Sigma_X=c$, that is, $X$ is stable.
\end{proof}

The above proposition will turn into a stronger result in the next section, when an additional hypothesis will be imposed. We close the section with a result which will be very useful later on. 

\begin{lemm}
\label{lemsrj}
Let $X=(A,B,I,J)\in\mathbf{B}$ and $A',B'\in {\rm End}(V')$ where $V'$ is a complex vector space. If $X$ is stable then
\begin{gather*}
\begin{matrix}
\phi: & {\rm Hom}(V',V)^{\oplus 2}\oplus {\rm Hom}(V',W) & \longrightarrow &{\rm Hom}(V',V)   \\
      & (a,b,j)  & \longmapsto     & Ab-bA'+aB'-Ba+Ij
\end{matrix}
\end{gather*}
is surjective.
\end{lemm}

\begin{proof}
The linear morphism $\phi$ is not surjective if and only if there is a nonzero $y\in{\rm Hom}(V,V')$ such that $y^\dagger\in(\im \phi)^\perp$. But $y^\dagger\in(\im \phi)^\perp$ if and only if, for every $(a,b,j)\in{\rm Hom}(V',V)^{\oplus 2}\oplus {\rm Hom}(V',W)$, holds
$$ 
{\rm Tr}(\phi(a,b,j)y)=0
$$
which is equivalent, for every $(a,b,j)\in{\rm Hom}(V',V)^{\oplus 2}\oplus {\rm Hom}(V',W)$, to the following equalities
\begin{align*}
{\rm Tr}((Ab-bA')y)={\rm Tr}((yA-A'y)b)&=0 \\
{\rm Tr}((aB'-Ba)y)={\rm Tr}(a(B'y-yB))&=0 \\
     {\rm Tr}((Ij)y)={\rm Tr}((yI)j)&=0 
\end{align*}
which hold if and only if 
\begin{equation}
\label{equabi}
yA-A'y=B'y-yB=yI=0.
\end{equation}
So if $\phi$ is not surjective, let $y\neq 0$ satisfying (\ref{equabi}). Then $yI=0$ yields $I(W)\subset\ker y$, while $yA=A'y$ implies that $\ker y$ is $A$-invariant and $y B=B'y$ implies $\ker y$ is also $B$-invariant. Since $\ker y\subsetneq V$, it follows that $X$ is not stable.
\end{proof}


\section{The ADHM Variety}

Now we introduce the main object of our study.

\begin{defi}
The set $\calv=\calv(r,c):=\mu^{-1}(0)$ is called the {\em  ADHM variety}, i.e., the variety of solutions of 
$$
[A,B] + IJ = 0
$$
which is called the {\em ADHM equation}.
\end{defi}

We keep the notation previously introduced, i.e., $\calv^{\rm st}$, $\calv^{\rm cs}$, $\calv^{\rm reg}$, $\calv^{\rm sj}$, $\calv^{\rm ts}$, $\calv^{(s)}$, $\calv^{[s]}$ are, respectively, the intersections of $\mathbf{B}^{\rm st}$, $\mathbf{B}^{\rm cs}$, $\mathbf{B}^{\rm reg}$, $\mathbf{B}^{\rm sj}$, $\mathbf{B}^{\rm ts}$, $\mathbf{B}^{(s)}$, $\mathbf{B}^{[s]}$ with $\calv$ for $0\leq s\leq c$. In particular, we have that $\calv^{(c)}=\calv^{\rm st}$, and, since $X=(A,B,I,J)\in\calv^{(0)}$ if and only if $I=[A,B]=0$,
$$ 
\calv^{[0]} = \calv^{(0)} = \mathcal{C}_c \times {\rm Hom}(V,W)
$$
where $\mathcal{C}_c$ is the variety of $c\times c$ commuting matrices. 

\begin{lemm}
\label{lemsis}
For any $X\in\calv$ hold
\begin{enumerate}
\item[(i)] $\Sigma_X=\im R(X)$; 
\item[(ii)] $R(X)$ is surjective iff $X$ is stable.
\end{enumerate}
\end{lemm}

\begin{proof}
We already know from Proposition \ref{l1} that $\Sigma_X\supseteq\im R(X)$; besides, (ii) immediately follows from (i). Therefore, we just have to prove that $\Sigma_X\subseteq\im R(X)$. To do that, it suffices to show that if $X=(A,B,I,J)$, then $\im R(X)$ is $A$- and $B$-invariant and $I(W)\subset\im R(X)$.

Clearly, $I(W)\subset\im R(X)$; let us now show that $\im R(X)$ is $A$ and $B$-invariant. Note that
\begin{align*}
A(\im R(X))&=\sum\limits_{0\leqslant k,l\leqslant c-1} \im A^{k+1}B^lI  \\
&= \sum\limits_{0\leqslant l\leqslant c-1} \im A^{c}B^lI + \sum_{\substack{1\leqslant k\leqslant c-1 \\ 0\leqslant l\leqslant c-1}} \im A^{k}B^lI. 
\end{align*}
The second factor is clearly within $\im R(X)$. For the first factor, use Cayley-Hamilton Theorem to express $A^{c}$ as a linear combination of lower powers of $A$; it then follows that this factor is also within $\im R(X)$. So $\im R(X)$ is $A$-invariant.

Using the ADHM equation, one can see that 
$$ 
BA^k = A^kB + \sum_{r+s=k-1}A^rIJA^s. 
$$
Therefore we have
\begin{align*}
B(\im R(X)) &=\sum\limits_{0\leqslant k,l\leqslant c-1} \im BA^{k}B^lI  \\
            &=\sum\limits_{0\leqslant k,l\leqslant c-1} \im A^kB^{l+1}I+
                  \sum_{\substack{0\leqslant k,l\leqslant c-1 \\ r+s=k-1}}\im A^rIJA^sB^lI \\
            &\subseteq \sum\limits_{0\leqslant k\leqslant c-1} \im A^kB^cI +
                           \sum_{\substack{0\leqslant k\leqslant c-1 \\ 1\leqslant l\leqslant c-1}} \im A^kB^lI+
                               \sum_{0\leqslant r\leqslant c-2}\im A^rI.
\end{align*}
The second and third factors are clearly within $\im R(X)$. For the first factor, use again Cayley-Hamilton Theorem to express this turn $B^c$ as lower powers of $B$ to conclude that this factor must lie within $\im R(X)$. So $\im R(X)$ is $B$-invariant as well and we are done.
\end{proof}

\begin{theo} 
\label{thmvar} 
The ADHM variety $\calv(r,c)$ is a set theoretic complete intersection which is irreducible if and only if $r\geq 2$. Moreover, the set $\calv(r,c)^{(s)}$ is an irreducible quasi-affine variety of dimension $2rc+c^2-(r-1)(c-s)$ which is nonsingular if and only if either $s=c$ or $s=c-1$. In particular, $\calv(r,c)^{\rm st}$ and $\calv(r,c)^{\rm cs}$ are nonsingular irreducible quasi-affine varieties of dimension $2rc+c^2$.
\end{theo}

\begin{proof}
Set $\mathbf{D}:={\rm Hom}(W^{\oplus c^2},V)$. For $0\le s\le c$, set also $P^{(s)}:=\{T\in\mathbf{D}~|~ {\rm rk}\, T=s \}$ and $P^{[s]}:=\{T\in\mathbf{D}~|~ {\rm rk}\, T\leq s \}$. Then $P^{(s)}$ is open within $P^{[s]}$ which is closed within $\mathbf{D}$. Now $R:\mathbf{B}\to\mathbf{D}$ is continuos and, from the prior lemma, $\calv^{(s)}=\calv\cap R^{-1}(P^{(s)})$ and $\calv^{[s]}=\calv\cap R^{-1}(P^{[s]})$. So $\calv^{(s)}$ is open within $\calv^{[s]}$ which is closed within $\calv$. Therefore $\calv^{(s)}$ is a quasi-affine variety. 

Since $\calv^{\rm st}=\calv^{(c)}$, it follows that $\calv^{\rm st}$ is a quasi-affine variety; from Proposition \ref{prpsrj} it is also nonsingular. Since $\calv^{\rm st}$ is a variety, the isomorphism $X\to X^{\star}$ in $\mathbf{B}$ restricts to an isomorphism from $\calv^{\rm st}$ onto its image. But its image is precisely $\calv^{\rm cs}$ because, first, $X^\star$ is costable iff $X$ is stable and, second, $\mu(X^\star)=0$ iff $\mu(X)=0$. The dimension of both varieties is $\dim\mathbf{B}-\dim{\rm End}(V)=2(c^2+rc)-c^2=2rc+c^2$.

Now we claim that $\calv(r,c)^{(s)}$ is the total space of a rank $(r+s)(c-s)$ bundle
\begin{equation}
\label{equfib}
\calv(r,c)^{(s)}\longrightarrow{\rm G}(s,c)\times\calv(r,s)^{\rm st}\times\mathcal{C}_{c-s}
\end{equation}
where $G$ is the Grassmannian and $\mathcal{C}$ the commuting matrices variety. In fact, if $X\in \calv(r,c)^{(s)}$ we may move $\Sigma_X$ to the vector space spanned by $e_1,\,e_2,\ldots,e_s$. So $X=(A,B,I,J)$ where
\begin{equation} \label{fcanonica}
A = \left( \begin{array}{cc} A_1 & A_2 \\ 0 & A_3 \end{array} \right) ~~~~ B = \left( \begin{array}{cc} B_1 & B_2 \\ 0 & B_3 \end{array} \right) ~~~~
 I = \left( \begin{array}{c} I_1 \\ 0 \end{array} \right) ~~~~
J = \left( \begin{array}{cc} J_1 & J_2 \end{array} \right) 
\end{equation}
and
\begin{align*}
[A_1,B_1]+I_1J_1&=0 \\
A_1B_2-B_1A_2+A_2B_3-B_2A_3+I_1J_2&=0 \\
[A_3,B_3]&=0
\end{align*}
with the requirement that $X|_{\Sigma_X}=(A_1,B_1,I_1,J_1)$ is stable. Now the first equation is the ADHM equation for $X|_{\Sigma_X}$, the third equation is the commuting matrices equation for $A_3$ and $B_3$. The movement of $\Sigma_X$ to the $s$-dimensional standard space is described by the Grassmannian $G(s,c)$ and we have a natural map as in (\ref{equfib}). The fiber $\calv(r,c)^{(s)}_P$ over $P=(\Sigma_X,X|_{\Sigma_X},(A_3,B_3))$ is the set of $(A_2,B_2,J_2)$ which satisfy the second equation. This is clearly an affine space. Since $X|_{\Sigma_X}$ is stable, it follows from Lemma \ref{lemsrj}  that 
\begin{align*}
\dim \calv(r,c)^{(s)}_P&=2s(c-s)+r(c-s)-s(c-s) \\
              &=(r+s)(c-s).
\end{align*}
This proves the claim. Therefore $\calv^{(s)}$ is nonsingular iff $\mathcal{C}_{c-s}$ is so, which happens iff either $s=c$ or $s=c-1$.

Now let us compute the dimension. We have
\begin{align*}
\dim \calv^{(s)}&=\dim {\rm G}(s,c)+\dim \calv(r,s)^{\rm st}+\dim \mathcal{C}_{c-s}+\dim \varphi_P \\
                &=s(c-s)+2rs+s^2+(c-s)+(c-s)^2+(r+s)(c-s) \\
                &=(c-s)(s+1+c+r)+2rs+s^2 \\
                &=rc+c^2+c-s+rs \\
                &=rc+(rc-rc)+c^2+c-s+rs \\
                &=2rc+c^2-(r-1)(c-s).
\end{align*}

Now $\calv$ is the union of the $\calv^{(s)}$, so its dimension is the highest among them, which is $2rc+c^2$. It follows that $\calv$ is a set theoretic complete intersection.

Finally, we will prove that $\calv(r,c)^{(s)}$ is irreducible for each $r$, $c$ and $0\le s\le c$. First, we argue that $\calv(r,c)^{(c)}=\calv(r,c)^{\rm st}$ is irreducible.

The case $r=1$ is rather special because if $X=(A,B,I,J)$ is a stable solution of the ADHM equation, then $J=0$ by \cite[Prp. 2.8(1)]{N}. It follows that
the closure of $\calv(1,c)^{\rm st}$ is given by $\mathcal{C}_{c}\times{\rm Hom}(W,V)$, which is clearly irreducible. Hence $\calv(1,c)^{\rm st}$ is also irreducible.  

On the other hand, if $c=1$, $A$ and $B$ are simply numbers, while $I$ and $J$ can be regarded as vectors in $\mathbb{C}^{r}$,
$$
I = (x_1,\dots,x_r)\ \ \ \ \ \ J = (y_1,\dots,y_r)
$$
and, in this way, the ADHM equation reduces to
$$
IJ = \sum_{i=1}^{r} x_iy_i = 0.
$$
We will show that $\calv(r,1)$ is irreducible if $r\geq 2$. In order to do this, assume $\sum_{i=1}^{r} x_iy_i=(a_0+\sum_i a_ix_i+b_iy_i)(a_0'+\sum_i a_i'x_i+b_i'y_i)$. Fix $i$. We have $a_ia_i'=0$. Assume, without loss in generality, that $a_i'=0$. Then $a_ia_j'=0$ for every $j$, $a_ib_j'=0$ for every $j\neq i$, and $a_ib_i'=1$. This implies all $a_j'$ vanish and all $b_j'$ with $j\neq i$ vanish as well. This is impossible comparing both sides of polynomial equality above unless $r=1$. So $\calv(r,1)=\mathbb{C}^2\times V$ where $V$ is an irreducible variety of $\mathbb{C}^{2r}$ if $r\geq 2$. Therefore $\calv(r,1)$ is irreducible if $r\geq 2$. Since $\calv(r,1)^{\rm st}$ is open within $\calv(r,1)$, it is irreducible as well if $r\geq 2$. But if $r=1$ we have already seen it is also irreducible.

For the case $r,c\geq 2$, write any $(a,b,i,j)\in\mathbf{B}$ as  
$$
a = \left( \begin{array}{cc} a_1 & a_2 \\ e & a_3 \end{array} \right) ~~~~ b = \left( \begin{array}{cc} b_1 & b_2 \\ f & b_3 \end{array} \right) ~~~~
 i = \left( \begin{array}{c} i_1 \\ g \end{array} \right) ~~~~
j = \left( \begin{array}{cc} j_1 & j_2 \end{array} \right) 
$$
where the sizes of the matrices involved are as in (\ref{fcanonica}). Set
$$
x=(a_1,b_1,i_1,j_1,a_3,b_3)\ \ \ \ \ y=(a_2,b_2,j_2)\ \ \ \ \ z=(e,f,g)
$$
and, keeping the notation of (\ref{fcanonica}), define
\begin{align*}
\phi_1(x)&=[A_1,b_1]+[a_1,B_1]+i_1J_1+I_1j_1 \\
\phi_2(x)&=a_1B_2-a_2B_2+i_1J_2+A_2b_3-B_2a_3 \\
\phi_3(y)&=A_1b_2-b_2A_3+a_2B_3-B_1a_2-I_1j_1 \\
\phi_4(z)&=A_2f-B_2e \\
\phi_5(z)&=eB_1-B_3e+A_3f-fA_1+gJ_1 \\
\phi_6(z)&=eB_2-fA_2+gJ_2
\end{align*}
So if $X$ is as in (\ref{fcanonica}) and $s=c-1$, then
\begin{gather*}
\begin{matrix}
D_X\mu : & \mathbf{B} & \longrightarrow &{\rm End}(V)   \\
      & (x,y,z)  & \longmapsto     & \left( \begin{array}{cc} \phi_1(x)+\phi_4(z) &  \phi_2(x)+\phi_3(y)\\ \phi_5(z) &  \phi_6(z)\end{array} \right)
\end{matrix}
\end{gather*}
Now, since $X|_{\Sigma_X}$ is stable, Proposition \ref{prpsrj} and Lemma \ref{lemsrj} imply, respectively, that $\phi_1$ and $\phi_3$ are surjective. So $D_X\mu$ is surjective iff 
\begin{gather*}
\begin{matrix}
\psi_X : & \mathbb{C}^{2c+r-2} & \longrightarrow & \mathbb{C}^{c}  \\
       & z  & \longmapsto     & (\phi_5(z),\phi_6(z)) 
\end{matrix}
\end{gather*}
is surjective. Set $Q:=\calv\setminus\calv^{\rm sj}$. Note that 
\begin{equation}
\label{equq12}
Q\cap\calv^{(c-1)}=V_1\cup V_2
\end{equation} 
where 
$$
V_1:=\{X\in\calv^{(c-1)}\,|\,\phi_5\ {\rm is\ nonsurjective}\} 
$$ 
$$
V_2:=\{X\in\calv^{(c-1)}\,|\,\phi_6\ {\rm is\ nonsurjective}\}.
$$
Let us study the codimensions of the above varieties within $\calv^{(c-1)}$ in order to get the one of $Q\cap\calv^{(c-1)}$. From (\ref{equfib}), we have that $\calv^{(c-1)}$ is the total space of a rank $r+c-1$ bundle on $T:=\mathbb{P}^{c-1}\times\calv(r,c-1)^{\rm st}\times\mathbb{C}^2$. Write $W:=\mathbb{P}^{c-1}\times\calv(r,c-1)^{\rm st}$. Now $V_1$ is a subvariety of
$$
V_3:=\left\{X\in\calv^{(c-1)}\ \bigg{|}\begin{array}{c}A_3\ {\rm is\ an\ eigenvalue\ of}\ A_1 \\ B_3\ {\rm is\ an\ eigenvalue\ of}\  B_1\end{array}\right\}
$$ 
and since $(A_2,B_2,J_2)$ do not appear in the relations which describe $V_3$, it follows that $V_3$ is a rank $r+c-1$ bundle over the subvariety of $T$ given by
$$
V_4:=\left\{((A_1,B_1,\ldots),(A_3,B_3))\in W\times \mathbb{C}^2\ \bigg{|}\begin{array}{c}A_3\ {\rm is\ an\ eigenvalue\ of}\ A_1 \\ B_3\ {\rm is\ an\ eigenvalue\ of}\  B_1\end{array}\right\}
$$
and so $\dim T-\dim V_4=2$, hence $\dim\calv^{(c-1)}-\dim V_3=2$ which implies that $\dim\calv^{(c-1)}-\dim V_1\geq 2$. On the other hand, 
$$
V_2=\{X\in\calv^{(c-1)}\,|\,A_2=B_2=J_2=0\}
$$
so $V_2\cong T$ and hence $\dim\calv^{(c-1)}-\dim V_2=r+c-1$. Thus, from (\ref{equq12}), we get that $\dim\calv^{(c-1)}-\dim (Q\cap\calv^{(c-1)})\geq 2$.
Therefore, computing the codimension of $Q$ within $\calv$ we have
\begin{align*}
\codim Q &=\min \{\codim (Q\cap\calv^{(s)})\}_{s=0}^{c-1} \\
         &=\min \{\codim (Q\cap\calv^{(c-1)}),\min\{\codim \calv^{(s)}\}_{s=0}^{c-2}\} \\
         &\geq  \min \{(r-1)+2,\min\{(r-1)(c-s)\}_{s=0}^{c-2}\} \\
         &=\min \{r+1,2(r-1)\}\geq 2.
\end{align*}

It follows from \cite[Thm. 1.3]{H1} that $\calv^{\rm sj}$ is the tangent cone of a conected variety, so $\calv^{\rm sj}$ is conected. Since it is nonsingular, it is also irreducible. So $\calv^{\rm st}$, being open within $\calv^{\rm sj}$, is irreducible as well. It follows that $\calv(r,c)^{(s)}$ is the total space of a vector bundle over an irreducible basis (given by ${\rm G}(s,c)\times\calv(r,s)^{\rm st}\times\mathcal{C}_{c-s}$), hence $\calv(r,c)^{(s)}$ is also irreducible (see for instance \cite[Lem. 2.8]{CTT}).

Furthermore, since $\calv$ is the closure of its nonsingular locus and this one is contained in $\calv^{\rm sj}$, it follows that $\calv$ is the closure of $\calv^{\rm sj}$ as well. But we have just seen that, if $r,c\geq 2$, then $\calv(r,c)^{\rm sj}$ is irreducible, so the same holds for $\calv(r,c)$. Besides, we have already got the irreducibility of $\calv(r,1)$ if $r\geq 2$. On the other hand, to check necessity, note that $\dim\calv(1,c)^{(0)}=2c+c^2=\dim\calv(1,c)^{\rm st}$. So the vaiety $\calv(1,c)$ has an open subset -- $\calv(1,c)^{\rm st}$ -- and a proper closed subset -- $\calv(1,c)^{(0)}$ -- both of the same dimension, which cannot happen unless it is reducible. We are finally done.
\end{proof}

\begin{rema}\rm
Since $\calv(r,c)$ is irreducible for $r\ge2$, we also conclude that $\calv^{\rm sj}(r,c)$ concides with the nonsingular locus of $\calv(r,c)$ in this case.
\end{rema}

\begin{rema}\rm
Note that $\calv(1,c)$ has at least two irreducible components, i.e., the closures of $\calv(1,c)^{\rm st}$ and $\calv(1,c)^{\rm cs}$ within $\calv(1,c)$, which, by \cite[Prp. 2.8(1)]{N}, are given by the conditions $J=0$ and $I=0$, respectively. It can easily be checked that $\calv(1,1)$ has precisely $2$ irreducible components and, using software tools, that $\calv(1,2)$ has exactly $3$ irreducible components. So it makes sense to ask if $\calv(1,c)$ has always $c+1$ irreducible components.
\end{rema}

\begin{rema}
\label{remstr}
{\rm Consider the sequence of inclusions
$$ (\calv^{\rm st}\cup\calv^{\rm cs}) \subset \calv^{\rm sj} \subset \calv^{\rm ts}. $$
We will show with an example with $r=c=2$ that they are, in general, strict ones. Consider the following ADHM datum $X=(A,B,I,J)$ given by
$$ 
A = \left( \begin{array}{cc} a_1 & a_2 \\ 0 & a_3 \end{array} \right) ~~~~ B = \left( \begin{array}{cc} b_1 & b_2 \\ 0 & b_3 \end{array} \right) ~~~~
 I = \left( \begin{array}{cc} i_1 & i_2 \\ 0 & 0 \end{array} \right) ~~~~
J = \left( \begin{array}{cc} 0 & j_2 \\ 0 & j_4 \end{array} \right) .
$$
It is easy to check that if $(a_1 - a_3)b_2 + (b_1 - b_3)a_2 + i_1j_2 + i_2j_4 = 0$, then $X$ satisfies the ADHM equation. Note also that $A(S),B(S),I(W)\subset S \subset \ker J$, where $S$ is the subspace generated by the vector $(1,0)$; thus $X$ is neither stable nor costable. 

The corresponding Jacobian matrix is given by
$$
D_X\mu=\left( \begin{array}{cc} M & N  \end{array} \right) 
$$
where
$$
M=
\left( \begin{array}{cccccccc} 
0 & 0       & -b_2     & 0  & 0 & 0       & a_2      & 0  \\
b_2& b_3-b_1   & 0       & -b_2&-a_2& a_1-a_3   & 0       & a_2  \\  
0 & 0       & b_1-b_3   & 0  & 0 & 0       & a_3-a_1   & 0  \\        
0 & 0       & b_2      & 0  & 0 & 0       & -a_2     & 0  
\end{array} \right) 
$$
$$
N=\left( \begin{array}{cccccccc} 
0 & 0 & 0     & 0  & i_1 & 0       & i_2      & 0      \\
j_2& j_4   & 0       & 0     & 0   & i_1       &0 &i_2    \\
0 & 0   & 0  & 0 & 0       & 0   & 0  & 0 \\
0 & j_2      & j_4  & 0 & 0       & 0     & 0 & 0
\end{array} \right) 
$$
Hence if in addition $a_1\ne a_3$, $b_1\ne b_3$, $i_1\ne0$ and $j_4\ne0$, then $X\in\calv^{\rm ns}$. This shows that the first inclusion is strict.

On the other hand, if $a_1=a_3$ and $b_1=b_3$, then $X\not\in\calv^{\rm ns}$. However, if $i_1,i_2\ne0$, $j_1,j_2\ne0$ and $a_1,b_1\ne0$, then $X\in\calv^{\rm ts}$.
This shows that the second inclusion is also strict, as claimed.}
\end{rema}


\subsection{Gieseker to Uhlenbeck map}

Now we define an analogue of the so-called \emph{Gieseker to Uhlenbeck map} (c.f. \cite[pp. 196-198]{B}) purely in terms of the ADHM data. 

\begin{prop}
\label{prpule}
There exists a natural map 
$$ \calv(r,c)^{\rm st}/GL(c)\longrightarrow \coprod_{s=0}^{c}\, \calv(r,s)^{\rm reg}/GL(s)\,\times\, \mathcal{C}_{c-s}/GL(c-s). $$
\end{prop}

Recall that $\calv(r,c)^{\rm st}/GL(c)$ (resp. $\calv(r,c)^{\rm reg}/GL(c)$) coincides with the moduli space of framed torsion free (resp. locally free) sheaves of rank $r$ and second Chern class $c$ on $\p2$. In our context, $\calv(r,c)^{\rm st}/GL(c)$ plays the role of the \emph{Gieseker compactification} of $\calv(r,c)^{\rm reg}/GL(c)$. On the other hand, $\calv(r,c)^{\rm reg}/GL(c)$ can also be interpreted as the moduli space of framed instantons on $\R^4$, and the space $\calv(r,s)^{\rm reg}/GL(s)\,\times\, \mathcal{C}_{c-s}/GL(c-s)$ is called the \emph{Uhlenbeck compactification} of $\calv(r,s)^{\rm reg}/GL(s)$ (recall that the quotient variety $\mathcal{C}_{c-s}/GL(c-s)$ coincides with the symmetric product $\text{Sym}^c(\A^2)$).

\begin{proof}
For any $X=(A,B,I,J)\in\calv(r,c)$, we define its costabilizing subspace as
$$
\Upsilon_X:=\langle S\subset V\,|\, A(S),B(S)\subset S\subset \ker J(S)\rangle,
$$
that is, the largest $A$- and $B$-invariant subspace of $V$ on which $J$ vanishes. We have that $X$ is costable if and only if
$\Upsilon_X=0$. 

Set $\calv(r,c)_{(s)}:=\{X\in\calv(r,c)\,|\,\dim\Upsilon_X=s\}$. If $X\in \calv(r,c)_{(s)}$ we may move $\Upsilon_X$ to the vector 
space spanned by $e_1,\,e_2,\ldots,e_s$. So $X=(A,B,I,J)$ where
$$
A = \left( \begin{array}{cc} A_1 & A_2 \\ 0 & A_3 \end{array} \right) ~~~~ B = \left( \begin{array}{cc} B_1 & B_2 \\ 0 & B_3 \end{array} \right) ~~~~
 I = \left( \begin{array}{c} I_1 \\ I_2 \end{array} \right) ~~~~
J = \left( \begin{array}{cc} 0 & J_2 \end{array} \right) 
$$
and
\begin{align*}
[A_1,B_1]&=0 \\
A_1B_2-B_1A_2+A_2B_3-B_2A_3+I_1J_2&=0 \\
[A_3,B_3]+I_2J_2&=0
\end{align*}
with $(A_3,B_3,I_2,J_2)$ costable. Thus $X':=(A_3,B_3,I_2,J_2)\in\calv(r,s)^{\rm cs}$ and also $Y:=(A_1,B_1)\in \mathcal{C}_{c-s}$. It is easily seen that if $X$ is stable, so is $X'$. Hence we are able to define the map $\overline{X}\mapsto (\overline{X'},\overline{Y})$.
\end{proof}


\subsection{The structure of the ADHM category}
\label{adhmcat}

\

\

In what follows, it will be important to put our subject within a categorical framework. As it is well-known, an ADHM datum $(A,B,I,J)$ can be regarded as representation of the so-called ADHM quiver
$$
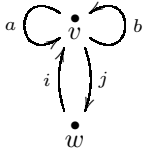
\xymatrix{
\stackrel{\bullet}{v} \ar@/^/[d]^{j} \ar@(ul,dl)[]_{a} \ar@(dr,ur)[]_{b} \\
\stackrel{\bullet}{w} \ar@/^/[u]^{i} } ~~~~ 
$$
provided with the relation $ab-ba+ij$. Let $\cala$ be the category of representations of the ADHM quiver. Recall that $\cala$ is an abelian category, and that for any representations $\mathbf{R}=(V,W,(A,B,I,J))$ and $\mathbf{R}'=(V',W',(A',B',I',J'))$ in $\cala$, we have that 
$$
{\rm Hom}_{\cala}(\mathbf{R},\mathbf{R}')=\left\{(f,g)\in{\rm Hom}(V,V')\oplus{\rm Hom}(W,W')\bigg{|}\begin{array}{ll}fA=A'f & fI=I'g \\ fB=B'f & gJ=J'f\end{array}\right\}.
$$
A representation $\mathbf{R}=(V,W,X)$ is said to be \emph{stable} (resp. \emph{costable}) if $X$ is stable (resp. costable) and we set
$$
\mathbf{S}_{\mathbf{R}}:=(\Sigma_X,W,X|_{\Sigma_X})
$$ 
to be the \emph{stable restriction} of $\mathbf{R}$. Set $N_X:=V/\Sigma_X$. If $X=(A,B,I,J)$ is as in (\ref{fcanonica}) we have that
\begin{equation}
\label{equzzr}
\mathbf{Z}_{\mathbf{R}}:=(N_X,\{0\},(A_3,B_3,0,0))
\end{equation} 
is the \emph{quotient representation} of $\mathbf{R}$. We also call the triple $(\dim W,\dim \Sigma_X, \dim N_X)$ the \emph{type vector} of $\mathbf{R}$.

Let us now pause for some general theory. Consider an additive category $\mathcal{C}$ and let $\mathcal{B}$ be a full subcategory of $\mathcal{C}$; recall that the {\em right orthogonal} to $\mathcal{B}$ is the full subcategory $\mathcal{B}^\perp$ of $\mathcal{C}$ consisting of the objects $D$ such that ${\rm Hom}_{\mathcal{C}}(B,D)=0$ for every $B$ in $\mathcal{B}$. 

\begin{defi}\label{radm}
Let $\mathcal{C}$ be an additive category. A pair of full subcategories $(\mathcal{B},\mathcal{D})$ is a called a {\em torsion pair} of $\mathcal{C}$ if $\mathcal{D}\subset\mathcal{B}^\perp$ and if for every object $C$ of $\mathcal{C}$, there exists a short exact sequence $0\to B \to C \to D\to 0$ with $B$ in $\mathcal{B}$ and $D$ in $\mathcal{D}$.
\end{defi}

Given an object $C$ of $\mathcal{C}$, the exact sequence $0\to B \to C \to D\to 0$ with $B$ in $\mathcal{B}$ and $D$ in $\mathcal{D}$ is unique up to isomorphism, so that 
the assignments $\iota^*(C):=B$ and $\jmath^*(C):=D$ yield additive functors $\iota^*:\mathcal{C}\to\mathcal{B}$ and $\jmath_*:\mathcal{C}\to\mathcal{D}$. It follows (see \cite[Proposition 1.2]{BR}) that $\iota^*$ is a right adjoint to the inclusion functor $\iota_*:\mathcal{B}\to\mathcal{C}$, while $\jmath^*$ is a left adjoint to the inclusion functor $\jmath_*:\mathcal{D}\to\mathcal{C}$. 

Let us now show how Definition \ref{radm} applies in our context. Let $\mathcal{S}$ be the full subcategory of $\cala$ whose objects are the stable representations, and let $\mathcal{Z}$ be the full subcategory of $\cala$ consisting of representations with type vector of the form $(0,s,l)$; note that any representation in $\mathcal{Z}$ is of the form $(V,\{0\},(A,B,0,0))$ where $A$ and $B$ are commuting operators on $\rm{End}(V)$.

\begin{prop}
\label{propcat}
The pair $(\mathcal{S},\mathcal{Z})$ is a torsion pair on $\mathcal{A}$.
\end{prop}

\begin{proof}
For any $\mathbf{R}$ one gets a short exact sequence 
$$
0\to \mathbf{S}_{\mathbf{R}} \to \mathbf{R} \to \mathbf{Z}_{\mathbf{R}} \to 0
$$ 
and clearly $\mathbf{S}_{\mathbf{R}}$ and $\mathbf{Z}_{\mathbf{R}}$ are in $\mathcal{S}$ and $\mathcal{Z}$, respectively.

Now we claim that $\mathcal{Z}\subset \mathcal{S}^\perp$. Take $\mathbf{S}=(V,W,X)\in \mathcal{S}$ with $X=(A,B,I,J)$ stable and let $\mathbf{Z}=(V',0,(A',B',0,0))$ be a representation in $\mathcal{Z}$. Then given $(f,g)\in{\rm Hom}_{\cala}(\mathbf{S},\mathbf{Z})$, we have that $fA=A'f$, $fB=B'f$ and $fI=0$. It follows that $\ker f$ is $A$ and $B$-invariant and it contains the image of $I$; since $X$ is stable, we have that $\ker f=V$ thus $f=0$; since $g\in{\rm Hom}(W,0)$, then also $g=0$, proving the claim.  
\end{proof}

As before, $\iota_*:\mathcal{S}\to\mathcal{A}$ and $\jmath_*:\mathcal{Z}\to\mathcal{A}$ be the inclusion functors, and let $\iota^*$ and $\jmath^*$ be their right and left adjoints, respectively. It this way, one has that $\mathbf{S}_{\mathbf{R}}=\iota^*(\mathbf{R})$ and $\mathbf{Z}_{\mathbf{R}}=\jmath^*(\mathbf{R})$.

\begin{rema}\rm
Recently, several authors have considered the so-called \emph{ADHM sheaves} (see e.g. \cite{Di,HL,Sz}), which are (twisted) representations of the ADHM quiver into ${\rm Coh}(Y)$, the category of coherent sheaves on a projective variety $Y$. These in turn may also be regarded as maps from $Y$ to the ADHM variety $\calv(r,c)$. One can show that if the image of $Y$ intersects $\calv(r,c)^{\rm st}$, then the corresponding ADHM sheaf is stable in the sense of \cite[Definition 2.2]{Di} and \cite[Definition 3.1]{HL}.
\end{rema}


\section{Perverse coherent sheaves on $\p2$}

In this Section we introduce, following Kashiwara \cite{Ka}, a $t$-structure on $D^{\rm b}(\p2)$, the bounded derived category of coherent sheaves on $\p2$, and define a functor from the ADHM category $\cala$ to the the core of this $t$-structure.

We start with notation. Let $Y$ be a finite dimensional, nonsingular, separated noetherian scheme. Let $\text{Mod}(\oy)$ be the category of $\oy$-modules and $D(\oy)$ be its derived category. Let also $D(Y)$ be the derived category of coherent sheaves on $Y$. We set $D_{\rm qc}(\oy)$ (resp. $D_{\rm coh}(\oy)$) to be the full triangulated subcategory of $D(\oy)$ consisting of complexes with quasi-coherent (resp. coherent) cohomology. Recall that $D^{\rm b}(Y)$ is naturally equivalent to $D_{\rm coh}^{\rm b}(\oy)$.

We will also use the costumary notation $D^{\le n}(\oy)$ to mean complexes $M^\bullet$ in $D(\oy)$ such that $H^p(M^\bullet)=0$ for all $p>n$; similarly, $D^{\ge n}(\oy)$ to mean complexes $M^\bullet$ in $D(\oy)$ such that $H^p(M^\bullet)=0$ for all $p<n$.

A \emph{family of supports} on $Y$ is a set $\Phi$ of closed subsets of $Y$ satisfying the following conditions: (i) if $Z\in\Phi$ and $Z'$ is a closed subset of $Z$, then $Z'\in\Phi$; (ii) if $Z,Z'\in\Phi$, then $Z\cup Z'\in\Phi$; (iii) $\emptyset\in\Phi$.

A \emph{support datum} on $Y$ is a decreasing sequence $\mathbf{\Phi}:=\{\Phi^n\}_{n\in\Z}$ of families of supports satisfying the following conditions: (i) for $n\ll0$, $\Phi^n$ is the set of all closed subsets of $Y$; (ii) $n\gg0$, $\Phi^n=\{\emptyset\}$.

Given a support datum on $Y$, Kashiwara introduces the following subcatgories of $D^{\rm b}_{\rm qc}(\oy)$:
$$ {}^{\mathbf{\Phi}}D_{\rm qc}^{\le n}(\oy) := \left\{ M^\bullet \in D^{\rm b}_{\rm qc}(\oy) ~|~ {\rm supp}(H^k(M^\bullet)) \in \Phi^{k-n} ~\forall k \right\}, $$
$$ {}^{\mathbf{\Phi}}D_{\rm qc}^{\ge n}(\oy) := \left\{ M^\bullet \in D^{\rm b}_{\rm qc}(\oy) ~|~ R\Gamma_{\Phi^k}(M^\bullet) \in D^{\ge k+n}(\oy) ~\forall k \right\}. $$
The functor $\Gamma_{\Phi}:{\rm Mod}(\oy) \to {\rm Mod}(\oy)$ is defined as follows:
$$ \Gamma_{\Phi} := \lim_{Z\in\Phi} \Gamma_Z(F). $$
Then one has, for each open subset $U\subset Y$:
\begin{equation}\label{gamma}
\Gamma_{\Phi}(F)(U) = \{ \sigma\in F(U) ~|~ \overline{{\rm supp}\sigma} \in \Phi \}.
\end{equation}

Finally, the \emph{support} (or \emph{perversity}) \emph{function} associated to the support datum $\mathbf{\Phi}$ (see \cite[Lem. 5.5]{Ka}) is:
\begin{gather*}
\begin{matrix}
p_{\mathbf{\Phi}} :& Y                             & \longrightarrow &  \Z \\
      & x & \longmapsto     & {\rm max}\{ n\in\Z ~|~ \overline{x}\in\Phi^n \}.
\end{matrix}
\end{gather*}

Now consider as in \cite[p. 857]{Ka}
$$
{}^{\mathbf{\Phi}}D_{\rm coh}^{\le 0}(\oy) := {}^{\mathbf{\Phi}}D_{\rm qc}^{\le 0}(\oy) \cap D^{\rm b}(Y)
$$
$$
{}^{\mathbf{\Phi}}D_{\rm coh}^{\ge 0}(\oy) := {}^{\mathbf{\Phi}}D_{\rm qc}^{\ge 0}(\oy) \cap D^{\rm b}(Y).
$$
It is shown in \cite[Theorem 5.9]{Ka} that if the support function $p_\Phi$ satisfies the following condition
\begin{equation}\label{pervcond}
p_{\mathbf{\Phi}}(y) - p_{\mathbf{\Phi}}(x) \le \codim(\overline{\{y\}}) - \codim(\overline{\{x\}})\ \text{if}\ y\in \overline{\{x\}},
\end{equation}
then $( {}^{\mathbf{\Phi}}D_{\rm coh}^{\le 0}(\oy),{}^{\mathbf{\Phi}}D_{\rm coh}^{\ge 0}(\oy) )$ defines a $t$-structure on $D^{\rm b}(Y)$. 

We now especialize to the case we are concerned in. Set $Y=\p2$, and fix, for the remainder of the paper, homogeneous coordinates $(x:y:z)$ in $\p2$; let
$$
\ell_\infty:=\{(x:y:z)\in\p2\,|\,z=0\};
$$
be the \emph{line at infinity}. Consider de following support datum $\mathbf{\Phi}=\{\Phi^k\}_{k\in\Z }$ with 
\begin{align*}
\Phi^k & := \{ \rm all~closed~subsets~of~\p2 \} ~~ {\rm for} ~ k\le0 \\
\Phi^1 & := \{ \rm all~closed~subsets~of~\p2\setminus\ell_\infty \} \\
\Phi^k & := \{ \emptyset \}  ~~ {\rm for} ~ k\ge2.
\end{align*}
The corresponding perversity function $p_\Phi : \p2 \to \Z$ is given by: $p_{\mathbf{\Phi}}(x)=0$ iff $\overline{x}\cap\ell_\infty\ne\emptyset$ and $p_{\mathbf{\Phi}}(x)=1$ otherwise; in other words, $p_{\mathbf{\Phi}}(x)=1$ iff $x$ is a closed point away from $\ell_\infty$, and $p_{\mathbf{\Phi}}(x)=0$ otherwise. One easily checks that such function does satisfy the condition (\ref{pervcond}).

\begin{defi}
A \emph{perverse coherent sheaf on} $\p2$ is an object in the core of the $t$-structure $( {}^{\mathbf{\Phi}}D_{\rm coh}^{\le 0}(\op),{}^{\mathbf{\Phi}}D_{\rm coh}^{\ge 0}(\op) )$ on $D^{\rm b}(\p2)$.
\end{defi}

\begin{lemm}\label{rcl}
If $M^\bullet$ is a perverse coherent sheaf on $\p2$, then:
\begin{itemize}
\item[(i)] $H^1(M^\bullet)$ is supported away from $\ell_\infty$;
\item[(ii)] $H^p(M^\bullet)=0$ for $p\ne0,1$;
\item[(iii)] $H^0(M^\bullet)$ has no sections $\sigma$ such that $\overline{{\rm supp}\sigma}$ does not intersect $\ell_\infty$.
\end{itemize}
\end{lemm}

\begin{proof}
If $M^\bullet \in {}^{\mathbf{\Phi}}D_{\rm coh}^{\le 0}(\op)$, then $H^p(M^\bullet)=0$ for $p\ge2$, and $H^1(M^\bullet)$ is supported outside $\ell_\infty$. If $M^\bullet \in {}^{\mathbf{\Phi}}D_{\rm coh}^{\ge 0}(\op)$, then $R\Gamma_{\Phi^0}(M^\bullet) \in D^{\ge 0}(\op)$. But $\Gamma_{\Phi^0}$ is just the identity functor on ${\rm Mod}(\op)$, thus $R\Gamma_{\Phi^0}(M^\bullet)=M^\bullet$, and $M^\bullet \in D^{\ge 0}(\op)$, that is, $H^p(M^\bullet)=0$ for $p\leq -1$.

To establish the last claim, we apply \cite[Lem. 3.3.(iii)]{Ka}: since $M^\bullet \in D^{\ge 0}(\op)$, it follows that
$H^0(R\Gamma_{\Phi^1}(M^\bullet))=\Gamma_{\Phi^1}(H^0(M^\bullet))$. But $R\Gamma_{\Phi^1}(M^\bullet)\in D^{\ge 1}(\op)$ since $M^\bullet \in {}^{\mathbf{\Phi}}D_{\rm coh}^{\ge 0}(\op)$, thus $\Gamma_{\Phi^1}(H^0(M^\bullet))=0$, and our claim follows from the definition of the functor $\Gamma_{\Phi^1}$, see (\ref{gamma}). 
\end{proof}

Notice that the fact that ${\rm supp}\ H^1(M^\bullet)$ does not intersect $\ell_\infty$ implies that it is $0$-dimensional, and hence finite; it thus makes sense to speak of the length of $H^1(M^\bullet)$, that is, $\dim H^0(H^1(M^\bullet))$.

\begin{defi}
The \emph{rank, charge} and \emph{length} of a perverse coherent sheaf $M^\bullet$ are defined as, respectively, $\rk (H^0(E^\bullet))$, $c_2(H^0(E^\bullet))$ and $\text{length}(H^1(E^\bullet))$.
\end{defi}

The (abelian) category of perverse coherent sheaves on $\p2$ will be denoted by $\calp$. We say that $M^\bullet \in\calp$ is \emph{trivial at infinity} if the restriction 
$H^0(M^\bullet)|_{\ell_\infty}$ is isomorphic to ${\mathcal O}_{\ell_\infty}^{\oplus r}$, where $r=\rk (H^0(M^\bullet))>0$; let $\calp_\infty$ denote the full subcategory of $\calp$ consisting of such objects. 

\begin{theo} \label{chperv}
$M^\bullet \in \calp_\infty$ if and only if
\begin{itemize}
\item[(i)] $H^p(M^\bullet)=0$ for $p\neq0,1$;
\item[(ii)] $H^0(M^\bullet)$ is a torsion free sheaf which is trivial at infinity;
\item[(iii)] $H^1(M^\bullet)$ is a torsion sheaf with support outside $\ell_\infty$.
\end{itemize}
\end{theo}

It follows from the above result that the definition of perverse coherent sheaves which are trivial at $\ell_\infty$ used in \cite[Section 5]{BFG} coincides with ours.

\begin{proof}
For the ``only if" part, by Lemma \ref{rcl}, it is enough to argue that $H^0(E^\bullet)$ is torsion free. Indeed, let $T$ be the torsion submodule of $H^0(E^\bullet)$. On one hand, the support of $T$ cannot intersect $\ell_\infty$, since $H^0(M^\bullet)|_{\ell_\infty}$ is locally free by hypothesis; on the other, item (iii) of Lemma \ref{rcl} implies that $T$ cannot be supported on points away from $\ell_\infty$ either. Thus $T=0$ and $H^0(E^\bullet)$ must be torsion free.

Conversely, let us first check that $M^\bullet\in{}^{\mathbf{\Phi}}D_{\rm coh}^{\le 0}(\op)$. This is quite clear, since ${\rm supp}(H^0(E^\bullet)) \in \Phi^{0}$, ${\rm supp}(H^1(E^\bullet)) \in \Phi^{1}$ and ${\rm supp}(H^k(M^\bullet))=\emptyset \in \Phi^{k}$ for $k\ge2$.

Second, let us check that $M^\bullet\in{}^{\mathbf{\Phi}}D_{\rm coh}^{\ge 0}(\op)$. Since $M^\bullet \in D_{\rm coh}^{\ge 0}(\op)$, then by \cite[Lem. 3.3.(iii)]{Ka} we conclude that $R\Gamma_{\Phi^k}(M^\bullet) \in D^{\ge 0}(\op)$ for every $k$. In particular, $R\Gamma_{\Phi^k}(M^\bullet) \in D^{\ge k}(\op)$ for every $k\le 0$. Since $R\Gamma_{\Phi^k}(M^\bullet)=0$ for every $k\ge2$ ($\Phi^k=\{\emptyset\}$ in this range), we also conclude that $R\Gamma_{\Phi^k}(M^\bullet) \in D^{\ge k}(\op)$ for every $k\ge 2$. Therefore, it only remains for us to show that $R\Gamma_{\Phi^1}(M^\bullet) \in D^{\ge 1}(\op)$, i.e., $H^0(R\Gamma_{\Phi^1}(M^\bullet))=0$. Again by \cite[Lem. 3.3(iii)]{Ka}, $H^0(R\Gamma_{\Phi^1}(M^\bullet))=\Gamma_{\Phi^1}(H^0(M^\bullet))$, so it suffices to argue that $\Gamma_{\Phi^1}(H^0(M^\bullet))=0$. But $H^0(M^\bullet)$ is a torsion free sheaf, so for every open set $U\subset\p2$ and every local section $\sigma\in H^0(M^\bullet)(U)$, we have that $\overline{{\rm supp}\sigma}=\p2$. Thus from (\ref{gamma}), we conclude that $\Gamma_{\Phi^1}(H^0(M^\bullet))=0$, as desired.
\end{proof}

Notice that rank, charge and length are the only topological invariants for objects in $\calp_\infty$, and if $r$, $s$, and $l$ are, respectively rank, charge and length of $M^\bullet\in\calp_\infty$, then ${\rm ch}(M^\bullet) = r - (s+l) h^2$.

\begin{rema}\rm 
In \cite[Section 2]{HL}, the authors introduce \emph{perverse instanton sheaves on $\p3$}, which are objects in the core of a t-structure on $D^{\rm b}(\p3)$ defined
through tilting on a torsion pair in ${\rm Coh}(\p3)$. A similar construction also applies to our case: our perverse coherent sheaves can also be regarded as objects in the core of a t-structure on $D^{\rm b}(\p2)$ defined through tilting on a torsion pair in ${\rm Coh}(\p2)$.
\end{rema}


\subsection{ADHM construction of perverse sheaves on $\p2$}
\label{secprv}

We will now construct a functor $\F:\cala\to\calp_\infty$ which extends the usual ADHM construction of instantons, as presented by Donaldson in \cite{D1} and later extended by Nakajima in \cite{N}; the link between representations of the ADHM quiver and perverse coherent sheaves was also discovered by Drinfeld, c.f. \cite[Thm. 5.7]{BFG}. We will further elaborated on Drinfeld's construction by discussing the role played by stabilizing subspaces, and also deriving relations between cohomologies of the complexes associated to a representation and its stable restriction and quotient representations.

\begin{defi}
\label{defcom}
\emph{Fix homogeneous coordinates $(x:y:z)$ in $\p2$. For any representation $\mathbf{R}=(V,W,(A,B,I,J))$ define the complex 
$$
E^\bullet_\mathbf{R} \,:~ 
V\otimes\op(-1) \stackrel{\alpha}{\longrightarrow} (V\oplus V\oplus W)\otimes\op \stackrel{\beta}{\longrightarrow} V\otimes\op(1)
$$
where
$$
\alpha = \left( \begin{array}{c} zA + x\id \\ zB + y\id \\ zJ \end{array} \right)\ \ \  
\beta = \left( \begin{array}{ccc}
-zB - y\id ~~ & ~~ zA + x\id ~~ & ~~ zI
\end{array} \right).
$$
Note that the ADHM equation is equivalent to $\beta\alpha=0$. Any complex on $\p2$ obtained in this way will be called an {\em ADHM complex}.}
\end{defi} 

Our aim now will be to show that any ADHM complex is a perverse coherent sheaf which is trivial at infinity, so the assignement $\mathbf{R} \mapsto E^\bullet_\mathbf{R}$ provides the desired functor. To see how the functor acts on morphisms, let the pair $(f,g)\in{\rm Hom}(V',V)\oplus{\rm Hom}(W',W)$ be a morphism between representations $\mathbf{R}'$ and $\mathbf{R}$; then one has the following morphism $\tilde{f}$ between the corresponding ADHM complexes $E^\bullet_{\mathbf{R}'}$ and $E^\bullet_\mathbf{R}$:
$$
\xymatrix{
V'\otimes\op(-1) \ar[r]^{\alpha'\ \ \ \ \ \ \ } \ar[d]^{f\otimes\id} & (V'\oplus V'\oplus W')\otimes\op \ar[r]^{\ \ \ \ \ \ 
\ \beta'} \ar[d]^{(f\oplus f\oplus g)\otimes\id} &
V'\otimes\op(1) \ar[d]^{f\otimes\id} \\
V\otimes\op(-1) \ar[r]^{\alpha\ \ \ \ \ \ \ } & (V\oplus V\oplus W)\otimes\op \ar[r]^{\ \ \ \ \ \ \ \beta} & V\otimes\op(1)
} . $$
So one defines $\F((f,g))$ to be the roof $E^\bullet_\mathbf{R}\stackrel{\id}{\leftarrow}E^\bullet_\mathbf{R}\stackrel{\tilde{f}}{\rightarrow}E^\bullet_{\mathbf{R}'}$.

Let us fix once and for all a representation $\mathbf{R}=(V,W,(A,B,I,J))$, and the corresponding ADHM complex $E^\bullet_\mathbf{R}$ as above.

\begin{lemm} 
\label{l2}
The sheaf map $\alpha$ is injective. The fiber maps $\alpha_P$ are injective for every $P\in\p2$ if and only if $\mathbf{R}$ is costable.
\end{lemm}

\begin{proof}
It is easily seen that $\alpha_P$ is injective for every $P\in\ell_\infty$. It follows that the set of points $P\in\p2$ for which $\alpha_P$ is not injective is a $0$-dimensional subscheme of $\p2$. The first claim of the lemma follows.

Now, if $\alpha_P$ is not injective for some $P=(p:q:1)\in\p2\setminus\ell_\infty$, then there is a nonzero vector $v\in V$ such that
$$ 
Av = - p v\ \ \ \ \ \ Bv = - q v\ \ \ \ \ \ Jv=0
$$
and hence $\langle v\rangle$ is a nonzero subspace of $V$ which is invariant under $A$, $B$ and contained in $\ker J$, so $\mathbf{R}$ is not costable.

Conversely, if $\mathbf{R}$ is not costable, let $S\subset V$ be a nonzero subspace satisfying $A(S),B(S)\subset S\subset \ker J$. The ADHM equation implies that $[A|_S,B|_S]=0$, so let $v\in S$ be a nonzero common $A$ and $B$ eigenvector with eigenvalues $p$ and $q$, respectively. It is then easy to see that if $P=(-p:-q:1)$ then $\alpha_P$ is not injective.
\end{proof}

The following is a well-known result (see \cite[Section 2.1]{N}), which we include here just for completeness.

\begin{lemm} \label{l3}
If $\mathbf{R}$ is stable, then $H^1(E^{\bullet}_{\mathbf{R}})=0$, and $H^0(E^{\bullet}_{\mathbf{R}})$ is a torsion free sheaf whose restriction to $\ell_\infty$ is trivial of rank $r=\dim W$ and second Chern class $c=\dim V$.
\end{lemm}

\begin{lemm}
\label{lemhh0}
For any $\mathbf{R}\in\cala$, holds $H^0(H^0(E^{\bullet}_{\mathbf{R}})(-1))=0$.
\end{lemm}

\begin{proof}
Set $E:= H^0(E^{\bullet}_{\mathbf{R}})$. Consider the two short exact sequences of sheaves:
\begin{equation}\label{um}
0 \to  V\otimes\op(-1) \stackrel{\alpha}{\longrightarrow} \ker\beta \to E \to 0
\end{equation}
\begin{equation}
\label{dois}
0 \longrightarrow \ker\beta \longrightarrow (V\oplus V\oplus W)\otimes\op \stackrel{\beta}{\longrightarrow} \im \beta \longrightarrow 0.
\end{equation}
From (\ref{um}), we get $H^0(E(-1))\cong H^0(\ker\beta(-1))$. On the other hand, from (\ref{dois}), we get $H^0(\ker\beta(-1))=0$.
\end{proof}

We are finally ready to establish the main result of this section.

\begin{theo} \label{mpt}
If $\mathbf{R}$ is a representation of the ADHM quiver with type vector $(r,s,l)$, then the associated ADHM complex $E^\bullet_{\mathbf{R}}$ is a perverse coherent sheaf on $\p2$ of rank $r$, charge $s$ and length $l$ which is trivial at infinity. Moreover, 
\begin{itemize}
\item[(i)] $H^0(E^\bullet_{\mathbf{R}})\simeq H^0(E^\bullet_{\mathbf{S}_{\mathbf{R}}})$;
\item[(ii)] $H^1(E^\bullet_{\mathbf{R}})\simeq H^1(E^\bullet_{\mathbf{Z}_{\mathbf{R}}})$;
\item[(iii)] $H^0(E^\bullet_{\mathbf{R}})$ is locally free if and only if $\mathbf{R}$ is costable.
\end{itemize}
\end{theo}

In categorical terms, one has that $H^0(\F(\mathbf{R}))\simeq\F(\iota^*(\mathbf{R}))$ and $H^1(\F(\mathbf{R}))\simeq\F(\jmath^*(\mathbf{R}))$, where $\iota^*$ and $\jmath^*$ are the adjoint functors introduced in the end of Section \ref{adhmcat}.

\begin{proof}
Lemma \ref{l2} implies $H^{-1}(E^\bullet_\mathbf{R})=0$, thus $H^{p}(E^\bullet_\mathbf{R})=0$ for $p\ne0,1$. For the remainder, set $\mathbf{S}:=\mathbf{S}_{\mathbf{R}}$, $\mathbf{Z}:=\mathbf{Z}_{\mathbf{R}}$, $\Sigma:=\Sigma_X$ and $N:=N_X$ to simply notation. One then has the following short exact sequence of complexes
\begin{equation}\label{diag}
\begin{array}{ccccc} 
 0 &  & 0 & & 0\\
\downarrow &  & \downarrow & & \downarrow \\
\Sigma\otimes\op(-1) & \stackrel{\alpha'}{\longrightarrow} & (\Sigma\oplus\Sigma\oplus W)\otimes\op & \stackrel{\beta'}{\longrightarrow} & \Sigma\otimes\op(1) \\
\downarrow &  & \downarrow & & \downarrow \\
V\otimes\op(-1) & \stackrel{\alpha}{\longrightarrow}  & (V\oplus V\oplus W)\otimes\op &  \stackrel{\beta}{\longrightarrow}  & V\otimes\op(1)  \\
\downarrow &  & \downarrow & & \downarrow \\
N\otimes\op(-1) &  \stackrel{\alpha''}{\longrightarrow}  & (N\oplus N)\otimes\op & \stackrel{\beta''}{\longrightarrow} & N\otimes\op(1) \\
\downarrow &  & \downarrow & & \downarrow \\
0  &  & 0  &  & 0 
\end{array}
\end{equation}
where, clearly, each line corresponds to $E^\bullet_{\mathbf{S}}$, $E^\bullet_{\mathbf{R}}$ and $E^\bullet_{\mathbf{Z}}$, respectively. 

Note that
$$
\alpha'' = \left( \begin{array}{c} zA_3 + x\id \\ zB_3 + y\id \end{array} \right).
$$
Thus $\alpha''$ is injective at $\ell_\infty$, hence it is injective as a sheaf map, which is equivalent to $H^{-1}(E^\bullet_{\mathbf{Z}})=0$. Moreover, $\beta'$ is surjective since $\mathbf{S}$ is stable (see Lemma \ref{l3} above) and so $H^1(E_{\mathbf{S}}^\bullet)=0$. Therefore, the long exact sequence of cohomology associated to the short exact sequence of complexes above simplifies to
$$ 
0 \to H^0(E^\bullet_{\mathbf{S}}) \to H^0(E^\bullet_{\mathbf{R}}) \to H^0(E^\bullet_\mathbf{Z}) \to 0 \to H^1(E^\bullet_{\mathbf{R}}) \to H^1(E^\bullet_\mathbf{Z})\to 0. 
$$
This implies that $H^1(E^\bullet_{\mathbf{R}}) \simeq H^1(E^\bullet_\mathbf{Z})$, and that $H^0(E^\bullet_{\mathbf{S}}) \simeq H^0(E^\bullet_{\mathbf{R}})$ if and only if $H^0(E^\bullet_\mathbf{Z})=0$. Now both $H^0(E^\bullet_\mathbf{Z})$ and $H^1(E^\bullet_\mathbf{Z})$ are torsion sheaves supported at finitely many points which are outside $\ell_{\infty}$. In fact, take $P=(x:y:0)\in\ell_{\infty}$ and write
$$
\alpha''_P = \left( \begin{array}{c} x\id \\ y\id \end{array} \right)\ \ \ \ \beta''_P=\left( \begin{array}{cc} -y\id & x\id \end{array} \right).
$$
Since $H^0(E^\bullet_\mathbf{Z})=\ker\beta''/\im\alpha''$ and $H^1(E^\bullet_\mathbf{Z})={\rm coker}\,\beta''$, the stalks of both sheaves vanish at $P$. Hence the supports of both sheaves are $0$-dimensional schemes because none of each meet $\ell_{\infty}$. In particular, $H^0(H^0(E^\bullet_\mathbf{Z}))=H^0(H^0(E^\bullet_\mathbf{Z})(-1))$ and the latter vanish from Lemma \ref{lemhh0}; thus $H^0(E^\bullet_\mathbf{Z})=0$, for it is supported at finitely many points. It follows that $H^0(E^\bullet_{\mathbf{S}}) \simeq H^0(E^\bullet_{\mathbf{R}})$, as desired. 

From Lemma \ref{l3} we conclude that $H^0(E^\bullet_{\mathbf{R}})$ is a torsion free sheaf which restricts trivially at infinity. Moreover, $H^1(E^\bullet_{\mathbf{R}})$ is a torsion sheaf supported away from the line $\ell_\infty$ since $H^1(E^\bullet_{\mathbf{Z}})$ does and we have just seen that they are isomorphic. It follows from Theorem \ref{chperv} that $E^\bullet_{\mathbf{R}}\in\calp_\infty$.

Lemma \ref{l3} also tell us that $E^\bullet_{\mathbf{R}}$ has rank $r$ and charge $s$. To see that its length is $l$, consider the exact sequences
\begin{equation}
\label{equnn1}
0 \longrightarrow \ker\beta'' \longrightarrow (N\oplus N)\otimes\op \stackrel{\beta''}{\longrightarrow} \im \beta'' \longrightarrow 0
\end{equation}
\begin{equation}
\label{equnn2}
0 \to  N\otimes\op(-1) \stackrel{\alpha''}{\longrightarrow} \ker\beta'' \to H^0(E^\bullet_{\mathbf{Z}}) \to 0
\end{equation}
\begin{equation}
\label{equnn3}
0 \longrightarrow \im \beta''\longrightarrow N\otimes\op(1) \longrightarrow H^1(E^\bullet_{\mathbf{Z}}) \longrightarrow 0.
\end{equation}
From (\ref{equnn1}) we get that $H^0(\im\beta''(-1))\cong H^1(\ker\beta''(-1))$ and, similarly, we also get $H^1(\im\beta''(-1))\cong H^2(\ker\beta''(-1))$. From (\ref{equnn2}) we get that $\ker\beta''\cong N\otimes\op(-1)$ because $H^0(E^\bullet_{\mathbf{Z}})=0$. So $H^0(\im\beta''(-1))=H^1(\im\beta''(-1))=0$. Thus, from (\ref{equnn3}), it follows that $N\otimes H^0(\op) \cong H^0(H^1(E^\bullet_{\mathbf{Z}})(-1))$. Hence ${\rm length}(H^1(E^\bullet_{\mathbf{Z}}))=\dim N$ and so ${\rm length}(H^1(E^\bullet_{\mathbf{R}}))=l$.

It remains for us to prove the last claim: we have $H^0(E^\bullet_{\mathbf{R}})=\ker\beta/\im\alpha$ is locally free if and only if the $\alpha_P$ are injective for all $P\in\p2$, which holds if and only if $\mathbf{R}$ is stable, owing to Lemma \ref{l2}. We are done.
\end{proof}

\begin{rema}\rm
One can also introduce the notion of \emph{framed perverse coherent sheaves}. A \emph{framing} on $E^\bullet\in\calp_\infty$ is a choice of trivialization of $H^0(E^\bullet)$, i.e. an isomorphism $\phi:H^0(E^\bullet)|_{\ell_\infty} \stackrel{\sim}{\rightarrow} \oo_{\ell_\infty}^{\oplus r}$. A \emph{framed perverse coherent sheaf} on $\p2$ is a pair $(E^\bullet,\phi)$ consisting of a perverse coherent sheaf $E^\bullet$ which is trivial at infinity, and a framing $\phi$ on $E^\bullet$. Two framed perverse coherent sheaves $(E^\bullet,\phi)$ and $(F^\bullet,\varphi)$ are isomorphic if there exists an isomorphism $\Psi:E^\bullet\to F^\bullet$ such that $\varphi=\phi\circ(H^0\Psi|_\ell)$, where $H^0\Psi$ denotes the induced map $H^0(E^\bullet)\to H^0(F^\bullet)$.

Now consider the functor from the category of schemes to the category of \linebreak groupoids, denoted $\mathbf{P}(r,c)$, that assigns to each scheme $S$ the groupoid whose objects are $S$-families of framed perverse coherent sheaves $E^\bullet$ on $\p2$ of rank $r$, charge $s$ and length $l$ such that $c=s+l$. Drinfeld has proved that such functor defines a stack isomorphic to the quotient stack $[\calv(r,c)/GL(c)]$, see \cite[Thm. 5.7]{BFG}. We therefore conclude that the moduli stack $\mathbf{P}(r,c)$ is irreducible if an only if $r\ge2$.

It follows that the quotient stack $[\calv(r,c)^{(s)}/GL(c)]$ is isomorphic to the functor $\mathbf{P}(r,c,l)$ that assigns to each scheme $S$ the groupoid whose objects are $S$-families of framed perverse coherent sheaves $E^\bullet$ on $\p2$ of rank $r$, charge $s$ and length $l$; such stacks are irreducible for each $r$, $s$ and $l$.

It is also worth mentioning that it follows from the proof of \cite[Thm. 5.7]{BFG} that the functor $\F:\cala\to\calp_\infty$ is essentially surjective.
\end{rema}

%
%
%



\begin{thebibliography}{99}

\bibitem{ADHM}
Atiyah, M., Drinfeld, V., Hitchin, N., Manin, Yu.:
Construction of instantons.
Phys. Lett. {\bf 65A}, 185-187 (1978)

\bibitem{B}
Baranovsky, V.:
{\em Moduli of Sheaves on Surfaces and Action of the Oscilator Algebra}.
J. Diff. Geom. {\bf 55} (2000), 193--227.

\bibitem{BR}
Beligiannis, A., Reiten, I.:
Homological and Homotopical aspects of torsion theories.
Mem. Am. Math. Soc. {\bf 188} (2007).

\bibitem{BFG}
Braverman, A., Finkelberg, M., Gaitsgory. D.: 
Uhlenbeck Spaces via Affine Lie Algebras.
Progr. Math. {\bf 244} (2004), 17--135.




\bibitem{CTT}
Coand\u{a}, I., Tikhomirov, A., Trautmann, G.:
{\em Irreducibility and smoothness of the moduli space of mathematical 5-instantons over $\mathbb{P}^3$}.
Int. J. Math. {\bf 14} (2003), 1--45.

\bibitem{Di}
Diaconescu, D.-E.:
Moduli of ADHM sheaves and local Donaldson-Thomas theory.
Preprint math/0801.0820.

\bibitem{D1}
Donaldson, S.:
Instantons and Geometric Invariant Theory.
Commun. Math. Phys. {\bf 93} (1984), 453--460.

\bibitem{FJ2}
Frenkel, I. B., Jardim, M.:
Complex ADHM equations, and sheaves on $\p3$.
J. Algebra {\bf 319} (2008), 2913--2937.


\bibitem{H1}
Hartshorne, R.:
Cohomological Dimension of Algebraic Varieties.
Ann. Math. {\bf 88} (1968), 403--450.

\bibitem{HL}
Hauzer, M., Langer, A.:
{\em Moduli Spaces of Framed Perverse Instanton Sheaves on} $\mathbb{P}^3$.
Preprint arXiv:1003.5859.

\bibitem{J}
Jardim, M.:
Atiyah--Drinfeld--Hitchin--Manin construction of framed instanton sheaves.
C. R. Acad. Sci. Paris, Ser. I {\bf 346} (2008), 427--430.

\bibitem{JVM}
Jardim, M., Martins, R. V.:
Linear and Steiner bundles on projective varieties.
To appear in Comm. Alg. (2010).

\bibitem{Ka}
Kashiwara, M.: 
$t$-structures on the derived categories of holonomic $\mathcal D$-modules
and coherent $\mathcal O$-modules.
Moscow Math. J. {\bf 4} (2004), 847--868.

\bibitem{KN}
Kronheimer, P. B., Nakajima, H.:
Yang-Mills instantons on ALE gravitational instantons.
Math. Ann. {\bf 288} (1990), 263--307.

\bibitem{N}
Nakajima, H.:
{\em Lectures on Hilbert schemes of points on surfaces}.
Providence: American Mathematical Society, 1999.

\bibitem{N1}
Nakajima, H.:
Instantons on ALE spaces, quiver varieties, and Kac-Moody algebras.
Duke Math. J. {\bf 76} (1994), 365--416.

\bibitem{NS}
Nekrasov, N., Schwarz, A.:
Instantons on noncommutative $\R^4$, and $(2,0)$ superconformal six-dimensional theory.
Comm. Math. Phys. {\bf 198} (1998), 689--703.

\bibitem{OSS}
Okonek, O., Schneider, M., Spindler, H.:
{\em Vector bundles on complex projective spaces}.
Boston: Birkhauser, 1980.

\bibitem{ST}
Siu, Y.-T., Trautmann, G.:
{\em Gap-sheaves and extension of coherent analytic subsheaves}.
Lec. Notes Math, {\bf 172}.
Berlin: Springer-Verlag (1971)

\bibitem{Sz}
Szendroi, B.:
Sheaves on fibered threefolds and quiver sheaves.
Comm. Math. Phys. {\bf 278} (2008), 627--641.

\bibitem{T}
Taniguchi, T.:
ADHM construction of super Yang-Mills instantons.
J. Geom. Phys. {\bf 59} (2009), 1199--1209.


\end{thebibliography}
\end{document}